\documentclass[leqno,11pt]{article}
\usepackage{amsfonts}
\usepackage{graphics,graphicx,amssymb}
\usepackage{amsmath,amsthm}
\numberwithin{equation}{section}
\newtheorem{theorem}{Theorem}[section]
\newtheorem{corollary}{Corollary}[section]
\newtheorem{proposition}{Proposition}[section]

\newtheorem{definition}{Definition}[section]

\newfont{\got}{eufm9 scaled 1095}

\newfont{\w}{msbm9 scaled\magstep1}

\begin{document}

\title{Linear Connections on Normal Almost Contact Manifolds with Norden Metric}

\footnotetext[1]{This work is partially supported by The Fund for
Scientific Research of the University of Plovdiv, Bulgaria, Project
RS09-FMI-003.}
\author{Marta Teofilova}

\date{}

\maketitle

\begin{abstract}
Families of linear connections are constructed on almost contact
manifolds with Norden metric. An analogous connection to the
symmetric Yano connection is obtained on a normal almost contact
manifold with Norden metric and closed structural 1-form. The
curvature properties of this connection are studied on two basic
classes of normal almost contact manifolds with Norden metric.

\end{abstract}

\section*{Introduction}

The geometry of the almost contact manifolds with Norden metric
($B$-metric) is a natural extension of the geometry of the almost
complex manifolds with Norden metric \cite{Gan} in the case of an
odd dimension. A classification of the almost contact manifolds with
Norden metric consisting of eleven basic classes is introduced in
\cite{Ga-Mi}.

An important problem in the geometry of the manifolds equipped with
an additional tensor structure and a compatible metric is the study
of linear connections preserving the structures of the manifold. One
such connection on an almost contact manifold with Norden metric,
namely the canonical connection, is considered in \cite{Al,Man-Gri}.

In the present work we construct families of linear connections on
an almost contact manifold with Norden metric, which preserve by
covariant differentiation some or all of the structural tensors of
the manifold. We obtain a symmetric connection on a normal almost
contact manifold with Norden metric, which can be considered as an
analogue to the well-known Yano connection \cite{Ya1,Ya2}.


\section{Preliminaries}

Let $(M,\varphi,\xi,\eta,g)$ be a $(2n+1)$-dimensional \emph{almost
contact manifold with Norden metric}, i.e. $(\varphi,\xi,\eta)$ is
an \emph{almost contact structure}: $\varphi$ is an endomorphism of
the tangent bundle of $M$, $\xi$ is a vector field, $\eta$ is its
dual 1-form, and $g$ is a pseudo-Riemannian metric, called a
\emph{Norden metric} (or a $B$-\emph{metric}), such that
\begin{equation}\label{1}
\begin{array}{l}
\varphi^2x = -x + \eta(x)\xi,\qquad \eta(\xi)=1,\medskip\\
g(\varphi x,\varphi y) = -g(x,y) + \eta(x)\eta(y)
\end{array}
\end{equation}
for all $x,y \in \mathfrak{X}(M)$.

From (\ref{1}) it follows $\varphi \xi=0$, $\eta \circ \varphi =0$,
$\eta(x)=g(x,\xi)$, $g(\varphi x,y)=g(x,\varphi y)$.

The associated metric $\tilde{g}$ of $g$ is defined by
$\tilde{g}(x,y)=g(x,\varphi y) + \eta(x)\eta(y)$ and is a Norden
metric, too. Both metrics are necessarily of signature $(n+1,n)$.

Further, $x,y,z,u$ will stand for arbitrary vector fields in
$\mathfrak{X}(M)$.

Let $\nabla$ be the Levi-Civita connection of $g$. The fundamental
tensor $F$ is defined by
\begin{equation}\label{F}
F(x,y,z) = g\left((\nabla_x \varphi )y,z\right)
\end{equation}
and has the following properties
\begin{equation}\label{Fp}
\begin{array}{l}
F(x,y,z)=F(x,z,y),\medskip\\
F(x,\varphi y,\varphi z)=F(x,y,z) -
F(x,\xi,z)\eta(y)-F(x,y,\xi)\eta(z).
\end{array}
\end{equation}
Equalities (\ref{Fp}) and $\varphi \xi=0$ imply $F(x,\xi,\xi)=0$.

Let $\{e_i,\xi\}$, ($i=1,2,...,2n$) be a basis of the tangent space
$T_pM$ at an arbitrary point $p$ of $M$, and $g^{ij}$ be the
components of the inverse matrix of $(g_{ij})$ with respect to
$\{e_i,\xi\}$.

The following 1-forms are associated with $F$:
\begin{equation*}
\theta(x)=g^{ij}F(e_i,e_j,x),\qquad
\theta^{\ast}(x)=g^{ij}F(e_i,\varphi e_j,x),\qquad
\omega(x)=F(\xi,\xi,x).
\end{equation*}

A classification of the almost contact manifolds with Norden metric
is introduced in \cite{Ga-Mi}. Eleven basic classes $\mathcal{F}_i$
($i=1,2,...,11$) are characterized there according to the properties
of $F$.


The Nijenhuis tensor $N$ of the almost contact structure
$(\varphi,\xi,\eta)$ is defined in \cite{Sa} by
$N(x,y)=[\varphi,\varphi](x,y)+d\eta(x,y)\xi$, i.e.
$N(x,y)=\varphi^2[x,y] + [\varphi x,\varphi y]- \varphi[\varphi
x,y]- \varphi [x,\varphi y] + (\nabla_x \eta)y.\xi - (\nabla_y
\eta)x.\xi$. The almost contact structure in said to be integrable
if $N=0$. In this case the almost contact manifold is called
\emph{normal} \cite{Sa}.

In terms of the covariant derivatives of $\varphi$ and $\eta$ the
tensor $N$ is expressed by $N(x,y) = (\nabla_{\varphi x}\varphi) y -
(\nabla_{\varphi y}\varphi )x - \varphi(\nabla_x \varphi)y +
\varphi(\nabla_y \varphi)x + (\nabla_x \eta) y. \xi - (\nabla_y
\eta) x. \xi$, where $(\nabla_x \eta)y=F(x,\varphi y,\xi)$. Then,
according to (\ref{F}), the corresponding tensor of type (0,3),
defined by $N(x,y,z)=g(N(x,y),z)$, has the form
\begin{equation}\label{N3}
\begin{array}{l}
N(x,y,z) = F(\varphi x,y,z) - F(\varphi y,x,z) - F(x,y,\varphi z)
\medskip\\
\phantom{N(x,y,z)} + F(y,x,\varphi z) + F(x,\varphi y,\xi)\eta(z) -
F(y,\varphi x,\xi)\eta(z).
\end{array}
\end{equation}

The condition $N=0$ and (\ref{N3}) imply
\begin{equation}\label{166}
F(x,y,\xi)=F(y,x,\xi),\qquad \omega=0.
\end{equation}

The 1-form $\eta$ is said to be closed if $\mathrm{d}\eta=0$, i.e.
if $(\nabla_x\eta)y=(\nabla_y\eta)x$. The class of the almost
contact manifolds with Norden metric satisfying the conditions $N=0$
and $\mathrm{d}\eta=0$ is
$\mathcal{F}_1\oplus\mathcal{F}_2\oplus\mathcal{F}_4\oplus\mathcal{F}_5\oplus\mathcal{F}_6$
\cite{Ga-Mi}.

Analogously to \cite{Gan}, we define an associated tensor
$\widetilde{N}$ of $N$ by
\begin{equation}\label{Na}
\begin{array}{l}
\widetilde{N}(x,y,z) = F(\varphi x,y,z) + F(\varphi y,x,z) -
F(x,y,\varphi z)
\medskip\\
\phantom{N^\ast(x,y,z)} - F(y,x,\varphi z) + F(x,\varphi
y,\xi)\eta(z) + F(y,\varphi x,\xi)\eta(z).
\end{array}
\end{equation}
From $\widetilde{N}=0$ it follows $F(\varphi x,\varphi
y,\xi)=F(y,x,\xi)$, $\omega=0$. The class with $\widetilde{N}=0$ is
$\mathcal{F}_3\oplus\mathcal{F}_7$ \cite{Ga-Mi}.

The curvature tensor $R$ of $\nabla$ is defined as usually by
\begin{equation*}
R(x,y)z=\nabla_x\nabla_y z - \nabla_y\nabla_x z - \nabla_{[x,y]}z,
\end{equation*}
and its corresponding tensor of type (0,4) is given by
$R(x,y,z,u)=g(R(x,y)z,u)$.

A tensor $L$ of type (0,4) is said to be curvature-like if it has
the properties of $R$, i.e. $L(x,y,z,u)=-L(y,x,z,u)=-L(x,y,u,z)$ and
$\underset{x,y,z}{\mathfrak{S}}L(x,y,z,u)=0$ (first Bianchi
identity), where $\mathfrak{S}$ is the cyclic sum by three
arguments.

A curvature-like tensor $L$ is called a $\varphi$-K\"ahler-type
tensor if $L(x,y,\varphi z,\varphi u)=-L(x,y,z,u)$.

\section{Connections on almost contact manifolds\\ with Norden metric}

Let $\nabla^\prime$ be a linear connection with deformation tensor
$Q$, i.e. $\nabla^\prime_x y = \nabla_x y + Q(x,y)$. If we denote
$Q(x,y,z)=g(Q(x,y),z)$, then
\begin{equation}\label{11}
g(\nabla^{\prime}_x y - \nabla_x y,z) = Q(x,y,z).
\end{equation}

\begin{definition}
\emph{A linear connection $\nabla^\prime$ on an almost contact
manifold is called an} almost $\varphi$-connection \emph{if
$\varphi$ is parallel with respect to $\nabla^\prime$, i.e. if
$\nabla^\prime \varphi = 0$.}
\end{definition}

Because of (\ref{F}), equality (\ref{11}) and $\nabla^\prime
\varphi=0$ imply the condition for $Q$
\begin{equation}\label{14}
F(x,y,z)=Q(x,y,\varphi z) - Q(x,\varphi y, z).
\end{equation}

\begin{theorem}\label{th1}
On an almost contact manifold with Norden metric there exists a
10-parametric family of almost $\varphi$-connections
$\nabla^{\prime}$ of the form (\ref{11}) with deformation tensor $Q$
given by

\begin{equation}\label{12}
\begin{array}{l}
Q(x,y,z) = \frac{1}{2}\{F(x,\varphi y,z) + F(x,\varphi
y,\xi)\eta(z)\} - F(x,\varphi z,\xi)\eta(y) \medskip\\
+t_1\{F(y,x,z)+F(\varphi y,\varphi x,z)- F(y,x,\xi)\eta(z)
-F(\varphi y,\varphi x,\xi)\eta(z)\medskip\\
-F(y,z,\xi)\eta(x) - F(\xi,x,z)\eta(y)+
\eta(x)\eta(y)\omega(z)\}\medskip\\
+t_2\{F(z,x,y)+F(\varphi z,\varphi x,y)-F(z,x,\xi)\eta(y)-F(\varphi
z,\varphi x,\xi)\eta(y)\medskip\\
-F(z,y,\xi)\eta(x)-F(\xi,x,y)\eta(z)+\eta(x)\eta(z)\omega(y)\}\medskip\\
+t_3\{F(y,\varphi x,z)-F(\varphi y,x,z)-F(y,\varphi
x,\xi)\eta(z)+F(\varphi y,x,\xi)\eta(z)\medskip\\
-F(y,\varphi z,\xi)\eta(x)-F(\xi,\varphi
x,z)\eta(y)+\eta(x)\eta(y)\omega(\varphi z)\}\medskip\\
+t_4\{F(z,\varphi x,y) - F(\varphi z,x,y)-F(z,\varphi
x,\xi)\eta(y)+F(\varphi z,x,\xi)\eta(y)\medskip\\
-F(z,\varphi y,\xi)\eta(x) - F(\xi,\varphi
x,y)\eta(z)+\eta(x)\eta(z)\omega(\varphi y)\}\medskip\\
+t_5\{F(\varphi y,z,\xi) + F(y,\varphi z,\xi)-\eta(y)\omega(\varphi
z)\}\eta(x)\medskip\\
+t_6\{F(\varphi z,y,\xi)+F(z,\varphi y,\xi) - \omega(\varphi
y)\eta(z)\}\eta(x)\medskip\\
+t_7\{F(\varphi y,\varphi
z,\xi)-F(y,z,\xi)+\eta(y)\omega(z)\}\eta(x)\medskip\\
+t_8\{F(\varphi z,\varphi y,\xi) -
F(z,y,\xi)+\omega(y)\eta(z)\}\eta(x)\medskip\\
+t_9\omega(x)\eta(y)\eta(z)+t_{10}\omega(\varphi
x)\eta(y)\eta(z),\qquad t_i\in\mathbb{R},\quad i=1,2,...,10.
\end{array}
\end{equation}
\begin{proof}
The proof of the statement follows from (\ref{14}), (\ref{12}) and
(\ref{Fp}) by direct verification that $\nabla^\prime \varphi=0$ for
all $t_i$.
\end{proof}
\end{theorem}

Let $N=\mathrm{d}\eta=0$. By (\ref{N3}), (\ref{166}) from Theorem
\ref{th1} we obtain

\begin{corollary} Let $(M,\varphi,\xi,\eta,g)$ be a
$\mathcal{F}_1\oplus\mathcal{F}_2\oplus\mathcal{F}_4\oplus\mathcal{F}_5\oplus\mathcal{F}_6$-manifold.
Then, the deformation tensor $Q$ of the almost $\varphi$-connections
$\nabla^\prime$ defined by (\ref{11}) and (\ref{12}) has the form
\begin{equation}\label{18}
\begin{array}{l}
Q(x,y,z) = \frac{1}{2}\{F(x,\varphi y,z) + F(x,\varphi
y,\xi)\eta(z)\} - F(x,\varphi z,\xi)\eta(y) \medskip\\
+s_1\{F(y,x,z)+F(\varphi y,\varphi x,z)\}
+s_2\{F(y,\varphi x,z)-F(\varphi y,x,z)\}\medskip\\
+s_3F(y,\varphi z,\xi) \eta(x) +s_4F(y,z,\xi)\eta(x),
\end{array}
\end{equation}
where $s_1=t_1+t_2$, $s_2=t_3+t_4$, $s_3=2(t_5+t_6)-t_3-t_4$,
$s_4=-t_1-t_2-2(t_7+t_8)$.
\end{corollary}

\begin{definition}
\emph{A linear connection $\nabla^\prime$ is said to be} almost
contact \emph{if the almost contact structure $(\varphi,\xi,\eta)$
is parallel with respect to it, i.e. if $\nabla^\prime
\varphi=\nabla^\prime \xi = \nabla^\prime \eta=0$.}
\end{definition}
Then, in addition to the condition (\ref{14}), for the deformation
tensor $Q$ of an almost contact connection given by (\ref{11}) we
also have
\begin{equation}\label{15}
F(x,\varphi y,\xi) = Q(x,y,\xi)=-Q(x,\xi,y).
\end{equation}

\begin{definition}
\emph{A linear connection on an almost contact manifold with Norden
metric $(M,\varphi,\xi,\eta,g)$ is said to be} natural \emph{if
$\nabla^\prime \varphi=\nabla^\prime \eta=\nabla^\prime g=0$.}
\end{definition}
The condition $\nabla^\prime g=0$ and (\ref{11}) yield
\begin{equation}\label{16}
Q(x,y,z)=-Q(x,z,y).
\end{equation}
From (\ref{Na}) and (\ref{12}) we compute
\begin{equation}\label{17}
\begin{array}{l}
Q(x,y,z) + Q(x,z,y)=\medskip\\
=(t_1+t_2)\{\widetilde{N}(y,z,\varphi x) -
\widetilde{N}(\xi,y,\varphi
x)\eta(z)-\widetilde{N}(\xi,z,\varphi x)\eta(y)\}\medskip\\
\hspace{0.04in}-(t_3+t_4)\{\widetilde{N}(y,z,x)-\widetilde{N}(\xi,y,x)\eta(z)-\widetilde{N}(\xi,z,x)\eta(y)\}\medskip\\
\hspace{0.04in}-(t_5+t_6)\{\widetilde{N}(\varphi^2z,y,\xi)+\widetilde{N}(z,\xi,\xi)\eta(y)\}\eta(x)\medskip\\
\hspace{0.04in}+(t_7+t_8)\{\widetilde{N}(\varphi
z,y,\xi)-\widetilde{N}(\varphi
z,\xi,\xi)\eta(y)\}\eta(x)\medskip\\
\hspace{0.04in}+2\{t_9\omega(x)+t_{10}\omega(\varphi
x)\}\eta(y)\eta(z).
\end{array}
\end{equation}

By (\ref{11}), (\ref{12}), (\ref{15}), (\ref{16}) and (\ref{17}) we
prove the following
\begin{proposition}\label{th2}
Let $(M,\varphi,\xi,\eta,g)$ be an almost contact manifold with
Norden metric, and let $\nabla^\prime$ be the 10-parametric family
of almost $\varphi$-connections defined by (\ref{11}) and
(\ref{12}). Then
\begin{description}

\item[(i)] $\nabla^\prime$ are almost contact iff
$t_1+t_2-t_9=t_3+t_4-t_{10}=0$;

\item[(ii)] $\nabla^\prime$ are natural iff
$t_1+t_2=t_3+t_4=t_5+t_6=t_7+t_8=t_9=t_{10}=0$.
\end{description}
\end{proposition}

Taking into account equation (\ref{N3}), Theorem \ref{th1},
Proposition \ref{th2}, and by putting $p_1=t_1=-t_2$,
$p_2=t_3=-t_4$, $p_3=t_5=-t_6$, $p_4=t_7=-t_8$, we obtain

\begin{theorem} On an almost contact manifold with Norden metric
there exists a 4-parametric family of natural connections
$\nabla^{\prime\prime}$ defined by
\begin{equation*}\label{nat}
\begin{array}{l}
g(\nabla^{\prime\prime}_x y - \nabla_x y,z)=\frac{1}{2}\{F(x,\varphi
y,z)+F(x,\varphi
y,\xi)\eta(z)\}-F(x,\varphi z,\xi)\eta(y)\medskip\\
\phantom{g(\nabla^{\prime\prime}_x y - \nabla_x
y,z)}+p_1\{N(y,z,\varphi x)+N(\xi,y,\varphi
x)\eta(z)+N(z,\xi,\varphi
x)\eta(y)\}\medskip\\
\phantom{g(\nabla^{\prime\prime}_x y - \nabla_x y,z)}+p_2\{N(z,y,x)+N(y,\xi,x)\eta(z)+N(\xi,z,x)\eta(y)\}\medskip\\
\phantom{g(\nabla^{\prime\prime}_x y - \nabla_x
y,z)}+p_3\{N(\varphi^2
z,y,\xi)+N(z,\xi,\xi)\eta(y)\}\eta(x)\medskip\\
\phantom{g(\nabla^{\prime\prime}_x y - \nabla_x
y,z)}+p_4\{N(y,\varphi z,\xi)+N(\varphi z,\xi,\xi)\eta(y)\}\eta(x).
\end{array}
\end{equation*}
\end{theorem}

Since $N=0$ on a normal almost contact manifold with Norden metric,
the family $\nabla^{\prime\prime}$ consists of a unique natural
connection on such manifolds
\begin{equation}\label{D}
\begin{array}{l}
\nabla^{\prime\prime}_x y = \nabla_x y +
\frac{1}{2}\{(\nabla_x\varphi)\varphi y + (\nabla_x \eta)y.\xi \} -
\nabla_x \xi.\eta(y),
\end{array}
\end{equation}
which is the well-known canonical connection \cite{Al}.

Because of Proposition \ref{th2}, (\ref{17}) and the condition
$\widetilde{N}=0$, the connections $\nabla^\prime$ given by
(\ref{11}) and (\ref{12}) are natural on a
$\mathcal{F}_3\oplus\mathcal{F}_7$-manifold iff $t_1=-t_2$ and
$t_3=-t_4$.

Let $(M,\varphi,\xi,\eta,g)$ be in the class
$\mathcal{F}_1\oplus\mathcal{F}_2\oplus\mathcal{F}_4\oplus\mathcal{F}_5\oplus\mathcal{F}_6$.
Then, $N=\mathrm{d}\eta=0$ and hence the torsion tensor $T$ of the
4-parametric family of almost $\varphi$-connections $\nabla^\prime$
defined by (\ref{11}) and (\ref{18}) has the form
\begin{equation}\label{19}
\begin{array}{l}
T(x,y,z) = s_1\left\{F(y,x,z) - F(x,y,z) + F(\varphi y,\varphi
x,z)\right.\medskip\\
\left.\phantom{T(x,y,z)}-F(\varphi x,\varphi y,z)\right\}
+\frac{1-4s_2}{2}\{F(x,\varphi y,z)-F(y,\varphi x,z)\}\medskip\\
\phantom{T(x,y,z)}- s_4\{F(x,z,\xi)\eta(y)-F(y,z,\xi)\eta(x)\}\medskip\\
\phantom{T(x,y,z)}-(1-s_2+s_3)\{F(x,\varphi z,\xi)\eta(y) -
F(y,\varphi z,\xi)\eta(x)\}.
\end{array}
\end{equation}
Then, from (\ref{19}) we derive that $T=0$ if and only if
$s_1=s_4=0$, $s_2=\frac{1}{4}$, $s_3=-\frac{3}{4}$. By this way we
prove
\begin{theorem}
On a
$\mathcal{F}_1\oplus\mathcal{F}_2\oplus\mathcal{F}_4\oplus\mathcal{F}_5\oplus\mathcal{F}_6$-manifold
there exists a symmetric almost $\varphi$-connection $\nabla^\prime$
defined by
\begin{equation}\label{20}
\begin{array}{l}
\nabla_x^\prime y = \nabla_xy
+\frac{1}{4}\left\{2(\nabla_x\varphi)\varphi y
+(\nabla_y\varphi)\varphi x -(\nabla_{\varphi y}\varphi)x
\right.\medskip\\
\left. \phantom{\nabla_x^*y}+ 2(\nabla_x\eta)y.\xi -
3\eta(x).\nabla_y \xi - 4\eta(y).\nabla_x \xi\right\}.
\end{array}
\end{equation}
\end{theorem}
Let us remark that the connection (\ref{20}) can be considered as an
analogous connection to the well-known Yano connection
\cite{Ya1,Ya2} on a normal almost contact manifold with Norden
metric and closed 1-form $\eta$.

\section{Connections on $\mathcal{F}_4\oplus\mathcal{F}_5$-manifolds}

In this section we study the curvature properties of the
4-parametric family of almost $\varphi$-connections $\nabla^\prime$
given by (\ref{11}) and (\ref{18}) on two of the basic classes
normal almost contact manifolds with Norden metric, namely the
classes $\mathcal{F}_4$ and $\mathcal{F}_5$. These classes are
defined in \cite{Ga-Mi} by the following characteristic conditions
for $F$, respectively:
\begin{equation}\label{F4}
\begin{array}{l}
\mathcal{F}_4: F(x,y,z) = -\frac{\theta(\xi)}{2n}\{g(\varphi
x,\varphi y)\eta(z) + g(\varphi x, \varphi
z)\eta(y)\},\qquad\hspace{0.03in}
\end{array}
\end{equation}
\begin{equation}\label{F5}
\begin{array}{l}
\mathcal{F}_5: F(x,y,z) = -\frac{\theta^{\ast}(\xi)}{2n}\{g(\varphi
x,y)\eta(z) + g(\varphi x, z)\eta(y)\}.\qquad\quad\hspace{0.03in}
\end{array}
\end{equation}

The subclasses of $\mathcal{F}_4$ and $\mathcal{F}_5$ with closed
1-form $\theta$ and $\theta^\ast$, respectively, are denoted by
$\mathcal{F}_4^0$ and $\mathcal{F}_5^0$. Then, it is easy to prove
that on a $\mathcal{F}_4^0\oplus\mathcal{F}_5^0$-manifold it is
valid:
\begin{equation}\label{22}
x\theta(\xi)=\xi\theta(\xi)\eta(x),\qquad
x\theta^\ast(\xi)=\xi\theta^\ast(\xi)\eta(x).
\end{equation}

Taking into consideration (\ref{F4}) and (\ref{F5}), from (\ref{11})
and (\ref{18}) we obtain
\begin{proposition}
Let $(M,\varphi,\xi,\eta,g)$ be a
$\mathcal{F}_4\oplus\mathcal{F}_5$-manifold. Then, the connections
$\nabla^\prime$ defined by (\ref{11}) and (\ref{18}) are given by
\begin{equation}\label{21}
\begin{array}{l}
\nabla_x^{\prime}y = \nabla_xy + \frac{\theta(\xi)}{2n}\{g(x,\varphi
y)\xi - \eta(y)\varphi x\} +\frac{\theta^{\ast}(\xi)}{2n}\{g(x,y)\xi
- \eta(y)x\}\medskip\\
\phantom{\nabla_x^{\prime}y}+\frac{\lambda\theta(\xi)+\mu\theta^\ast(\xi)}{2n}\{\eta(x)y-\eta(x)\eta(y)\xi\}+
\frac{\mu\theta(\xi)-\lambda\theta^\ast(\xi)}{2n}\eta(x)\varphi y,
\end{array}
\end{equation}
where $\lambda=s_1+s_4$, $\mu=s_3-s_2$.
\end{proposition}
The Yano-type connection (\ref{20}) is obtained from (\ref{21}) for
$\lambda=0$, $\mu=-1$.

Let us denote by $R^\prime$ the curvature tensor of $\nabla^\prime$,
i.e. $R^\prime(x,y)z = \nabla^\prime_x\nabla^\prime_y z -
\nabla^\prime_y\nabla^\prime_x z - \nabla^\prime_{[x,y]}z$. The
corresponding tensor  of type (0,4) with respect to $g$ is defined
by $R^\prime(x,y,z,u)=g(R^\prime(x,y)z,u)$. Then, it is valid
\begin{proposition}
Let $(M,\varphi,\xi,\eta,g)$ be a
$\mathcal{F}_4^0\oplus\mathcal{F}_5^0$-manifold, and $\nabla^\prime$
be 2-parametric family of almost $\varphi$-connections defined by
(\ref{21}). Then, the curvature tensor $R^\prime$ of an arbitrary
connection in the family (\ref{21}) is of $\varphi$-K\"ahler-type
and has the form
\begin{equation}\label{25}
\begin{array}{l}
R^{\prime} = R +\frac{\xi\theta(\xi)}{2n}\pi_5 +
\frac{\xi\theta^{\ast}(\xi)}{2n}\pi_4+
\frac{\theta(\xi)^2}{4n^2}\{\pi_2 - \pi_4\}\medskip\\
\phantom{R^\prime} + \frac{\theta^{\ast}(\xi)^2}{4n^2}\pi_1 -
\frac{\theta(\xi)\theta^\ast(\xi)}{4n^2}\{\pi_3 - \pi_5\},
\end{array}
\end{equation}
where the curvature-like tensors $\pi_i$ ($i=1,2,3,4,5$) are defined
by \emph{\cite{Man-Gri}:}
\begin{equation}\label{26}
\begin{array}{l}
\pi_1(x,y,z,u) = g(y,z)g(x,u) - g(x,z)g(y,u),\medskip\\
\pi_2(x,y,z,u) = g(y,\varphi z)g(x,\varphi u) - g(x,\varphi
z)g(y,\varphi u),\medskip\\
\pi_3(x,y,z,u) = -g(y,z)g(x,\varphi u) + g(x,z)g(y,\varphi
u)\medskip\\
\phantom{\pi_3(x,y,z,u)} - g(x,u)g(y,\varphi z) + g(y,u)g(x,\varphi
z),\medskip\\
\pi_4(x,y,z,u)=g(y,z)\eta(x)\eta(u) - g(x,z)\eta(y)\eta(u)
\medskip\\
\phantom{\pi_4(x,y,z,u)}+ g(x,u)\eta(y)\eta(z) -
g(y,u)\eta(x)\eta(z),\medskip\\
\pi_5(x,y,z,u) = g(y,\varphi z)\eta(x)\eta(u) - g(x,\varphi
z)\eta(y)\eta(u)\medskip\\
\phantom{\pi_5(x,y,z,u)}+g(x,\varphi u)\eta(y)\eta(z) - g(y,\varphi
u)\eta(x)\eta(z).
\end{array}
\end{equation}
\begin{proof}
It is known that the curvature tensors of two linear connections
related by an equation of the form (\ref{11}) satisfy
\begin{equation}\label{24}
\begin{array}{l}
g(R^\prime(x,y)z,u) = R(x,y,z,u) + (\nabla_x Q)(y,z,u) - (\nabla_y
Q)(x,z,u)
\medskip\\
\phantom{g(R^\prime(x,y)z,u)}+ Q(x,Q(y,z),u) - Q(y,Q(x,z),u).
\end{array}
\end{equation}
Then, (\ref{25}) follows from  (\ref{24}), (\ref{21}), (\ref{22})
and (\ref{26}) by straightforward computation.
\end{proof}
\end{proposition}

Let us remark that (\ref{25}) is obtained in \cite{Man-Gri} for the
curvature tensor of the canonical connection (\ref{D}) on a
$\mathcal{F}_4^0\oplus\mathcal{F}_5^0$-manifold, i.e. the connection
(\ref{21}) for $\lambda=\mu=0$.

\smallskip

\begin{tabbing}
  \emph{Marta Teofilova}\\
  \emph{University of Plovdiv}\\
  \emph{Faculty of Mathematics and Informatics}\\
  \emph{236 Bulgaria Blvd.}\\
  \emph{4003 Plovdiv, Bulgaria}\\
  \texttt{e-mail:\ marta@uni-plovdiv.bg}\\
  \end{tabbing}

\end{document}